\documentclass[a4paper,12pt]{article}
\usepackage[francais,english]{babel}
\usepackage{latexsym,amssymb}
\usepackage{hyperref}
\usepackage{enumerate}
\usepackage{amsmath}
\usepackage{pb-diagram}
\newtheorem{thm}{Theorem}[section]
\newtheorem{cor}[thm]{Corollary}

\newtheorem{definition}[thm]{Definition}
\newtheorem{pro}[thm]{Proposition}
\newtheorem{eexemples}[thm]{Examples}

\def\spin{{\rm{Spin}}}

\def\SU{{\rm{SU}}}

\newcommand{\HH}{\mathbb{H}}
\newcommand{\C}{\mathbb{C}}
\newcommand{\Z}{\mathbb{Z}}
\newcommand{\R}{\mathbb{R}}
\newcommand{\T}{\mathbb{T}}
\newcommand{\s}{\mathbb{S}}
\newcommand{\ad}{\mathrm{Ad}}

\def \be{\begin{eqnarray*}}
\def \ee{\end{eqnarray*}}
\def \ben{\begin{enumerate}}
\def \een{\end{enumerate}}
\def \beit{\begin{itemize}}
\def \eeit{\end{itemize}}
\def \bui#1#2{\mathrel{\mathop{\kern 0pt#1}\limits^{#2}}}
\def \buil#1#2{\mathrel{\mathop{\kern 0pt#1}\limits_{#2}}}
\def \bflr{\begin{flushright}}
\def \eflr{\end{flushright}}

\def \lmt{\longmapsto}
\def \ovl{\overline}

\newcommand{\pa}[1]{\left(#1\right)}

\begin{document}
\title{{\bf A spectral estimate for the Dirac operator on Riemannian f\mbox{}lows}}
\author{{\small{\bf Nicolas Ginoux}}\footnote{Fakult\"at f\"ur Mathematik,
Universit\"at Regensburg,
D-93040 Regensburg,
E-mail: nicolas.ginoux@mathematik.uni-regensburg.de}, {\small\bf{Georges Habib}}\footnote{Lebanese University, Faculty of sciences II, Mathematics Department, P.O. Box 90656 Fanar-Matn, Lebanon, 
E-mail: ghabib@ul.edu.lb}}
\date{\today}
\maketitle

\begin{abstract}
\noindent
We give a new upper bound for the smallest eigenvalues of the Dirac operator on a Riemannian flow carrying transversal Killing spinors. We derive an estimate on Sasakian and on $3$-dimensional manifolds and partially classify those satisfying the limiting case. 
Finally, we compare our estimate with a lower bound in terms of a natural tensor depending on the eigenspinor.
\end{abstract}

\noindent\begin{small}{\it Mathematics Subject Classif\mbox{}ication}: 53C12, 53C25, 53C27, 58J50, 35P15\\
\end{small}

\noindent\begin{small}{\it Key words}: Foliations, Sasakian manifolds, Spin Geometry, Spectral geometry, Estimation of eigenvalues - upper and lower bounds\\
\end{small}
$ $\\

\section{Introduction}
The spectrum of the Dirac operator has been studied for a long time. Lower bounds of the Dirac operator are in general obtained by a suitable modification of the Levi-Civita connection and the use of the Schr\"odinger-Lichnerowicz formula. The limiting cases are characterized by the existence of special spinors which give rise to particular geometries.\\

\noindent For a manifold $M^n$ isometrically immersed into one of the three simply-connected space-forms $\R^{n+1}, \s^{n+1}, \HH^{n+1},$ C. B\"ar got upper bounds for the eigenvalues of the Dirac operator of $M$ in terms of the mean curvature of the hypersurface \cite{Baer98}. His results follow from the min-max-principle using parallel or Killing spinors as test-spinors.\\

\noindent In this paper, we aim at studying the spectrum of the Dirac operator on manifolds arising as total space of submersions over real space forms. More generally, we study {\it Riemannian flows} (see Section \ref{sec:pre}) which are locally given by Riemannian submersions with $1$-dimensional fibres. 
We are interested in the following question:\\\\ 
\noindent{\it Can one relate the spectrum of the Dirac operator on a spin manifold submerged onto the space-forms to geometric quantities?}\\\\  
\noindent The most natural situation to start with consists in assuming the Riemannian flow $(M,\mathcal{F})$ to carry what we call an $(\alpha,\beta)$-transversal Killing spinor and which can be thought of as the lift of some Killing spinor on the base, see \cite{GinHabib1} for an account of geometrical properties of flows with transversal Killing spinors. Using this spinor field as a test-spinor, we derive eigenvalue estimates for the fundamental Dirac operator of a closed manifold $M$ in terms of the {\it O'Neill tensor} \cite{ONeill66} of the flow, which is the natural geometric tensor expected in this context. We begin with the general framework in Section \ref{sec:mainthm}, then we focus our attention on the particular cases where $M$ is Sasakian or $3$-dimensional. If $M$ is Sasakian (Section \ref{sec:Sasaki}), then the general estimate can be substantially simplified and provides the existence of harmonic spinors for suitable deformations of the metric. 
In Section \ref{pro:dimentrois}, we give a complete classification of $3$-dimensional manifolds satisfying the limiting case: we show that the O'Neill tensor is constant and hence the manifold is either a local Riemannian product or homothetic to a Sasakian manifold.\\

\noindent Section \ref{sec:comp} sheds a new light on the estimate due to O. Hijazi in terms of the energy-momentum tensor \cite{Hij95}. Indeed we show that, if the limiting case of our estimate is attained on $3$-dimensional flows, then so is Hijazi's lower bound (see Proposition \ref{estimdime3}) whereas this fact is still true on Sasakian manifolds for special sections of the spinor bundle of $M$ (see Proposition \ref{tcomparHijazi}). To illustrate this, we treat a lot of important examples.\\\\

\noindent{\bf Acknowledgement.} The authors thank the universities of Nancy and Potsdam, the Max-Planck Institute for Mathematics in the Sciences as well as the Sonderforschungsbereich 647 ``{\it Raum - Zeit - Materie. Analytische und Geometrische Strukturen}'' of the Deutsche Forschungsgemeinschaft for their support in the preparation of \cite{GinHabib1} and this paper. It is a pleasure to thank Bernd Ammann, Christian B\"ar and Oussama Hijazi for valuable comments.
We also thank the referee for his/her speedy but pertinent review.

\section{Preliminaries} \label{sec:pre}
\setcounter{equation}{0}
\noindent
For preliminaries about spin structures on Riemannian foliations, we refer to \cite{Habibthese}.
Let $(M^{n+1},g,\mathcal{F})$ be an $(n+1)$-dimensional Riemannian manifold together with a Riemannian flow $\mathcal{F}$ defined by a unit vector field $\xi$ \cite{Carriere84}. It is a $1$-dimensional foliation of $M$ satisfying 
\begin{equation}
(\mathcal{L}_\xi g)(Z,W)=0,
\label{eq:12}
\end{equation}
for all $Z,W$ orthogonal to $\xi$, where $\mathcal{L}_\xi$ is the Lie derivative in the direction of $\xi.$ The metric $g$ is a bundle-like metric in the sense of \cite{Reinh59}. We denote by $Q:=\xi^\perp$ the normal bundle with the induced metric $g$.
The condition (\ref{eq:12}) gives rise to a natural covariant derivative $\nabla$ on $Q$, called the {\it transversal Levi-Civita connection} \cite{Tondeurlivre} and which is defined for any $Z\in\Gamma(Q)$ by
\begin{equation*}
\nabla _{X} Z :=
\left\{\begin{array}{ll}
\pi [X,Z], &  \textrm {$X=\xi$},\\\\
\pi (\nabla_{X}^{M}Z), & \textrm {$X\perp \xi$},
\end{array}\right.
\end{equation*}
where $\nabla^{M}$ is the Levi-Civita connection of $M$ and $\pi:TM\longrightarrow Q$ is the orthogonal projection. 
The connection $\nabla^M$ is then related to $\nabla$ through the following Gauss-type formula \cite{Habibthese}: for all $Z,W\in \Gamma(Q)$,
\begin{equation}
\left\{\begin{array}{ll}
\nabla^M_Z W=\nabla_Z W-g(h(Z),W)\xi, &\textrm {}\\\\
\nabla^M_\xi Z=\nabla_\xi Z+h(Z)-\kappa(Z)\xi,&\textrm {}
\end{array}\right.
\label{eq:22}
\end{equation}
where the tensor $h(Z):=\nabla^M_Z \xi$ is the O'Neill tensor and $\kappa:=\nabla^M_\xi \xi$ is the mean curvature of the flow. We point out that $h$ is a skew-symmetric endomorphism-field of $Q$ as a consequence of \eqref{eq:12}. In particular, one may associate a $2$-form $\Omega$ to $h$ on $Q$ through $\Omega(Z,W):=g(h(Z),W)$ for all $Z,W\in \Gamma(Q).$\\\

\noindent Local Riemannian products of $1$-dimensional with Riemannian manifolds constitute the simplest examples of Riemannian flows. In that case the unitary vector field $\partial_t$ is parallel, in particular it defines a Riemannian flow with geodesic fibres (i.e., $\kappa=0$) and vanishing O'Neill tensor. The transversal Levi-Civita connection locally corresponds to the one on the $1$-codimensional leaf. Other particular examples are provided by Sasakian ma\-ni\-folds, whose definition is recalled:
 
\begin{definition}\label{def:06}
A \emph{Sasakian} manifold is a Riemannian f\mbox{}low $(M^{n+1},g,\mathcal{F})$ given by a unit vector f\mbox{}ield $\xi$ with
\beit\item[i)] $\kappa=0$, i.e., the f\mbox{}low is minimal,
\item[ii)] $h^2=-\mathrm{Id}_Q$, i.e., $h$ is an almost-Hermitian structure on $Q$,
\item[iii)] $\nabla h=0$, i.e., $h$ is parallel on $Q$ (hence is a K\"ahler structure on $Q$).
\eeit
\end{definition}


\noindent It can be easily checked that this def\mbox{}inition is equivalent to the usual one, where one requires $\xi$ to be a unit Killing vector f\mbox{}ield satisfying
\[
\left|\begin{array}{ll}(\nabla^M\xi)^2&=-\mathrm{Id}_{TM}+\xi^\flat\otimes \xi\\
(\nabla_X^M\nabla^M\xi)(Y)&=g(\xi,Y)X-g(X,Y)\xi\end{array}\right.
\]
for all $X,Y\in \Gamma(TM)$. From Def\mbox{}inition \ref{def:06} the tensor field  $h$ defines a canonical K\"ahler structure on the normal bundle $Q$ of any Sasakian manifold. In particular such a manifold is always odd-dimensional. We shall from now on denote $m:=\frac{n}{2}$. In the following we shall also omit to write $\mathcal{F}$ for the f\mbox{}low and consider a Sasakian manifold as a triple $(M^{2m+1},g,\xi)$.

   

\begin{definition}\label{defetaEinstein}
A Riemannian flow is called \emph{$\eta$-Einstein} if and only if there exist real constants $\lambda$ and $\nu$ such that
\[\mathrm{Ric}_M=\lambda g+\nu\xi^\flat\otimes\xi^\flat.\]
\end{definition}

\noindent These metrics were first introduced by Okumura \cite{Okumura62} and were also studied by Galicki and Boyer \cite{BoyerGalicki00,BoyerGalicki06}. The scalar curvature of $\eta$-Einstein Sasakian manifolds is constant and equal to $2m(\lambda+1).$\\

\noindent The relation \eqref{eq:22} between the Levi-Civita connections can be easily extended to sections of $\Sigma M$: for every $\psi\in \Gamma(\Sigma M),$ we have \cite{Habibthese}
 \begin{equation}\label{eq:25}\left\{\begin{array}{ll}\nabla^M_\xi\psi=\nabla_\xi\psi+\frac{1}{2}\Omega\cdot\psi+\frac{1}{2}\xi\cdot\kappa\cdot \psi, &\textrm {}\\\\\nabla^M_Z\psi=\nabla_Z\psi+\frac{1}{2}\xi\cdot h(Z)\cdot\psi, &\textrm {}\end{array}\right.
\end{equation}
where $Z\in\Gamma(Q)$.
Here we may view \eqref{eq:25} as an analogue to the standard Gauss-Weingarten formula for hypersurfaces \cite{Trautman92}, where the O'Neill tensor plays the role of the second fundamental tensor.\\
On a Sasakian manifold, the form $\Omega$ is the K\"ahler form of the normal bundle and its action induces the following orthogonal decomposition of the spinor bundle of $M$, (see e.g. \cite{FriedKim00} for details): 
\begin{equation}\label{eqdecSigmaSasaki}
\Sigma M=\oplus_{r=0}^m \Sigma_r M,
\end{equation}
where $\Sigma_r M$ is the rank-$\begin{pmatrix}
m\\r
\end{pmatrix}$-eigenbundle associated with the eigenvalue $i(2r-m)$ of $\Omega$. Moreover, the action of $\xi$ on each $\Sigma_r M$, $0\leq r\leq m$, is given by 
\[\xi\cdot\psi_r=(-1)^{r+1}i\psi_r\]
for all $\psi_r\in\Gamma(\Sigma_rM)$. In particular, for all $Z\in\Gamma(Q)$, for all $\psi_0\in \Gamma(\Sigma_0 M)$ and $\psi_m\in \Gamma(\Sigma_m M),$  we have 
$$h(Z)\cdot \psi_0=i Z\cdot\psi_0\quad\text{and}\quad h(Z)\cdot \psi_m=-i Z\cdot \psi_m.$$


\section{Main theorem}\label{sec:mainthm}
\setcounter{equation}{0}
\noindent
In this section, we establish an upper bound for the first eigenvalue of the Dirac operator on a Riemannian flow carrying transversal Killing spinors (see \eqref{eqkilling} below) by computing the associated Rayleigh quotient. We first recall the Min-Max principle \cite{Chavellivre} for the Dirac operator on a compact Riemannian spin manifold $(M,g)$. We denote by $(\cdot\,,\cdot):=\int_M\Re\langle\cdot\,,\cdot\rangle v_g$ the scalar pro\-duct on $L^2(\Sigma M)$ induced by the Hermitian product $\langle\cdot\,,\cdot\rangle$ on $\Sigma M$ and by $\|\cdot\|$ its associated norm.\\

\noindent{\bf Min-Max Principle}: {\it Let $(\lambda_k)_{k\geq 1}$ be the spectrum of the Dirac operator on $(M,g)$ with $0\leq |\lambda_1|\leq\ldots\leq|\lambda_k|\leq|\lambda_{k+1}|\leq\ldots$. For any natural integer $k\geq 1,$ we have
$$\lambda_k^2=\min\limits_{E_k\subset \Gamma(\Sigma M)} \left\{\max\limits_{\psi\in E_k\setminus\{0\}} \{\frac{(D_M^2\psi,\psi)}{(\psi,\psi)}\}\right\},$$
where the minimum is taken on all $k$-dimensional vector subspaces $E_k$ of $\Gamma(\Sigma M).$}\\
Applying this theorem means choosing a subspace $E_k$ of sections of $\Sigma M$ called {\it test-sections}, on which the Rayleigh quotient $\frac{(D_M^2\psi,\psi)}{(\psi,\psi)}$ is evaluated.\\

\noindent Let $(M^{n+1},g,\mathcal{F})$ be a spin Riemannian flow. Recall \cite{GinHabib1} that, for $\alpha,\beta\in \C$, an $(\alpha,\beta)$-transversal Killing spinor on $M$ is a smooth section $\psi\in \Gamma(\Sigma M)$ satisfying, for all $Z\in \Gamma(Q)$,
\begin{equation}\label{eqkilling}
\left\{\begin{array}{ll}
\nabla_{\xi}\psi&=\alpha\,\xi\cdot\psi,\\
 &\\
\nabla_{Z}\psi&=\beta\, \xi\cdot Z\cdot\psi.
\end{array}\right.
\end{equation}

\noindent An $(\alpha,\beta)$-transversal Killing spinor is a parallel section of $\Sigma M$ w.r.t. the covariant derivative
\[X\lmt\nabla_X\varphi-\alpha g(X,\xi)\xi\cdot\varphi-\beta\xi\cdot X\cdot\varphi-\beta g(X,\xi)\varphi\]
for all $X\in\Gamma(TM)$ and $\varphi\in\Gamma(\Sigma M)$, in particular if it vanishes at one point it vanishes everywhere on $M$.
General geometric aspects of Riemannian flows carrying transversal Killing spinors have been studied in \cite{GinHabib1}, where in particular the existence of a large family of examples is shown, including circle bundles over manifolds with Killing spinors or suitable deformations of Einstein-Sasaki manifolds. Note also that a $(0,\beta)$-transversal Killing spinor corresponds to a basic Killing spinor \cite{Habibthese}.\\ 

\begin{thm}\label{lem:2}
Let $(M^{n+1},g,\mathcal{F})$ be a compact spin Riemannian flow. Assume the existence of a non-zero $(\alpha,\beta)$-transversal Killing spinor $\psi$ on $M$ for complex constants $\alpha,\beta$. Then
\begin{eqnarray}\label{eq:DM2psi}
\nonumber D_M^2\psi&=&(\alpha^2+n^2\beta^2+\frac{|\kappa|^2}{4})\psi-\frac{1}{4}\Omega\cdot\Omega\cdot\psi+\alpha\xi\cdot\Omega\cdot\psi+2\beta\Omega\cdot\psi\\
\nonumber& &+\frac{1}{4}\sum_{j,k=1}^n\xi\cdot e_j\cdot e_k\cdot\nabla_{e_j}h(e_k)\cdot\psi+\frac{1}{2}\nabla^M_{\xi}\Omega\cdot\psi-2n\alpha\beta\xi\cdot\psi\\
& &-\frac{1}{2}d\kappa^\flat\cdot\psi+\xi\cdot h(\kappa)\cdot\psi-\frac{1}{2}\mathrm{div}_M(\kappa)\psi+\beta\xi\cdot \kappa\cdot\psi.
\end{eqnarray}

\noindent In particular, if $\alpha$ and $\beta$ are real and $|\psi|=1$ on $M$, then we have an upper bound for the lowest eigenvalue $\lambda_1(D_M^2)$ of the square of the Dirac operator: 
\begin{eqnarray}\label{estimation}
\lambda_1(D_M^2)&\leq&\alpha^2+n^2\beta^2+\frac{\int_M|\kappa|^2v_g}{4\mathrm{Vol}(M)}+\frac{\|\Omega\cdot\psi\|^2}{4\mathrm{Vol}(M)}
+\alpha\frac{\int_M\Re\langle\xi\cdot\Omega\cdot\psi,\psi\rangle v_g}{\mathrm{Vol}(M)}\nonumber\\
&&+\frac{1}{2}\int_M\Re \langle\xi\cdot\kappa\cdot\Omega\cdot\psi,\psi\rangle v_g.
\end{eqnarray}

\end{thm}

\noindent Here it is important to notice that the upper bound (\ref{estimation}) depends on the $(\alpha,\beta)$-transversal Killing spinor $\psi$, which belongs to the data provided by the geometry of the flow. This should not be confused with Hijazi's eigenvalue bound (\ref{eq:0034}) depending on an eigenspinor; in general, $(\alpha,\beta)$-transversal Killing spinors are not eigenvectors of the Dirac operator.\\

\noindent{\bf Proof.} Plugging Equations \eqref{eqkilling} into \eqref{eq:25}, we compute the Dirac operator of $\psi$ and we obtain
\begin{eqnarray}\label{eqDMpsi}
\nonumber D_M\psi&=&\xi\cdot\nabla_\xi^M\psi+\sum_{i=1}^n e_j\cdot\nabla_{e_j}^M\psi\\
&=&-\alpha\psi+n\beta\xi\cdot\psi-\frac{1}{2}\kappa\cdot\psi-\frac{1}{2}\xi\cdot\Omega\cdot\psi.
\end{eqnarray}
We now compute $D_M^2\psi.$ On the one hand, we have
\be
\nabla_\xi^M(D_M\psi)&=&-\alpha\nabla_\xi^M\psi+n\beta\nabla_\xi^M\xi\cdot\psi+n\beta\xi\cdot\nabla_\xi^M\psi\\
& &-\frac{1}{2}\pa{(\nabla_\xi^M\kappa)\cdot\psi+\kappa\cdot\nabla_\xi^M\psi}\\
& &-\frac{1}{2}\pa{(\nabla_\xi^M\xi)\cdot\Omega\cdot\psi+\xi\cdot\nabla_\xi^M\Omega\cdot\psi+\xi\cdot\Omega\cdot\nabla_\xi^M\psi}\\
&=&-(\alpha^2+\frac{|\kappa|^2}{4})\xi\cdot\psi-n\alpha\beta\psi+\frac{n\beta}{2}\xi\cdot\Omega\cdot\psi+\frac{n\beta}{2}\kappa\cdot\psi\\
& &-\frac{1}{2}(\nabla_\xi^M\kappa)\cdot\psi+\frac{1}{2}h(\kappa)\cdot\psi-\frac{1}{2}\kappa\cdot\Omega\cdot\psi\\
& &-\frac{1}{2}\xi\cdot\nabla_\xi^M\Omega\cdot\psi-\frac{1}{4}\xi\cdot\Omega\cdot\Omega\cdot\psi.
\ee
On the other hand, we write for all $Z\in \Gamma(Q)$
\be
\nabla_Z^M(D_M\psi)&=&-\alpha\nabla_Z^M\psi+n\beta\nabla_Z^M\xi\cdot\psi+n\beta\xi\cdot\nabla_Z^M\psi\\
& &-\frac{1}{2}\pa{(\nabla_Z^M\kappa)\cdot\psi+\kappa\cdot\nabla_Z^M\psi}\\
& &-\frac{1}{2}\pa{(\nabla_Z^M\xi)\cdot\Omega\cdot\psi+\xi\cdot\nabla_Z^M\Omega\cdot\psi+\xi\cdot\Omega\cdot\nabla_Z^M\psi}\\
&=&-\alpha\beta\xi\cdot Z\cdot\psi-n\beta^2 Z\cdot\psi-\frac{\alpha}{2}\xi\cdot h(Z)\cdot\psi+\frac{n\beta}{2}h(Z)\cdot\psi\\
& &-\frac{1}{2}(\nabla_Z^M\kappa)\cdot\psi+\frac{\beta}{2}\xi\cdot\kappa\cdot Z\cdot\psi+\frac{1}{4}\xi\cdot\kappa\cdot h(Z)\cdot\psi\\
& &-\frac{1}{2}h(Z)\cdot\Omega\cdot\psi-\frac{1}{2}\xi\cdot\nabla_Z^M\Omega\cdot\psi+\frac{\beta}{2}\Omega\cdot Z\cdot\psi\\
& &+\frac{1}{4}\Omega\cdot h(Z)\cdot\psi.
\ee
Hence, we deduce that
\be
D_M^2\psi&=&\xi\cdot\nabla_\xi^M(D_M\psi)+\sum_{j=1}^ne_j\cdot\nabla_{e_j}^M(D_M\psi)\\
&=&(\alpha^2+n^2\beta^2+\frac{|\kappa|^2}{4}+\frac{|h|^2}{2})\psi-2n\alpha\beta\xi\cdot\psi+2\beta\Omega\cdot\psi+\beta\xi\cdot\kappa\cdot\psi\\
& &-\frac{1}{2}\xi\cdot(\nabla_\xi^M\kappa)\cdot\psi-\frac{1}{2}\sum_{j=1}^ne_j\cdot (\nabla_{e_j}^M\kappa)\cdot\psi+\frac{1}{4}\xi\cdot(\Omega\cdot\kappa-\kappa\cdot\Omega)\cdot\psi\\
& &+\frac{1}{2}\nabla_\xi^M\Omega\cdot\psi+\frac{1}{2}\sum_{j=1}^n\xi\cdot e_j\cdot\nabla_{e_j}^M\Omega\cdot\psi-\frac{1}{4}\Omega\cdot\Omega\cdot\psi\\
& &+\alpha\xi\cdot\Omega\cdot\psi+\frac{1}{2}\xi\cdot h(\kappa)\cdot\psi.
\ee
Using the local expression of the form $\Omega,$ it can be easily proved that 
$$\nabla_Z^M\Omega\cdot\psi=h^2(Z)\cdot\xi\cdot\psi+\frac{1}{2}\sum_{j=1}^ne_j\cdot\nabla_Zh(e_j)\cdot\psi.$$
Finally, we obtain
\begin{eqnarray}\label{eqDM2psinablah}
\nonumber D_M^2\psi\nonumber&=&(\alpha^2+n^2\beta^2+\frac{|\kappa|^2}{4})\psi-2n\alpha\beta\xi\cdot\psi+2\beta\Omega\cdot\psi+\beta\xi\cdot\kappa\cdot\psi\\
\nonumber& &-\frac{1}{2}d\kappa^\flat\cdot\psi+\frac{1}{2}\mathrm{div}_M(\kappa)\psi+\xi\cdot h(\kappa)\cdot\psi+\frac{1}{2}\nabla_\xi^M\Omega\cdot\psi\\
&&+\frac{1}{4}\sum_{j,k=1}^n \xi\cdot e_j\cdot e_k\cdot\nabla_{e_j}h(e_k)\cdot\psi-\frac{1}{4}\Omega\cdot\Omega\cdot\psi\nonumber\\
& &+\alpha\xi\cdot\Omega\cdot\psi.
\end{eqnarray}
From now on, we suppose that $\alpha$ and $\beta$ are real. Since the norm of $\psi$ is constant, we may assume that $|\psi|=1.$ We evaluate the Rayleigh quotient on $\psi$. Using the fact that the product $\langle X\cdot\psi,\psi\rangle$ is purely imaginary and that $\Re\langle X\cdot\psi,Y\cdot\psi\rangle=g(X,Y)|\psi|^2,$ for all $X,Y \in \Gamma(TM),$ the scalar product of \eqref{eqDM2psinablah} with $\psi$ gives
\be
\Re\langle D_M^2\psi,\psi\rangle &=&\alpha^2+n^2\beta^2+\frac{|\kappa|^2}{4}+\frac{1}{2}\mathrm{div}_M(\kappa)\\
& &+\frac{1}{4}\sum_{j,k=1}^n \Re\langle\xi\cdot e_j\cdot e_k\cdot\nabla_{e_j}h(e_k)\cdot\psi,\psi\rangle +\frac{1}{4}|\Omega\cdot\psi|^2\\
& &+\alpha\Re\langle\xi\cdot\Omega\cdot\psi,\psi\rangle.
\ee
Now it is a straightforward computation to show that
\[\frac{1}{4}\sum_{j,k=1}^n \Re \pa{\langle\xi\cdot\! e_j\cdot\! e_k\cdot\!\nabla_{e_j}h(e_k)\cdot\psi,\psi\rangle}=-\frac{1}{2}\mathrm{div}_M(W_\psi)+\frac{1}{2}\Re \pa{\langle \xi\cdot\!\kappa\cdot\!\Omega\cdot\!\psi,\psi\rangle},\]
where $W_\psi\in\Gamma(TM)$ is def\mbox{}ined by $g(W_\psi,X):=\Re \pa{\langle \xi\cdot\Omega\cdot\psi,X\cdot\psi\rangle}$ and in our convention $\mathrm{div}_M(X):=-\mathrm{tr}(\nabla^MX)$ for all $X\in \Gamma(TM)$.
Integrating on $M$ and applying the divergence theorem we f\mbox{}ind  
\be
\frac{\pa{D_M^2\psi,\psi}}{\pa{\psi,\psi}}&=&\alpha^2+n^2\beta^2+\frac{\int_M|\kappa|^2v_g}{4\mathrm{Vol}(M)}+\frac{\int_M\Re(\langle\xi\cdot\kappa\cdot\Omega\cdot\psi,\psi\rangle) v_g}{2\mathrm{Vol}(M)}\\
&&+\frac{\|\Omega\cdot\psi\|^2}{4\mathrm{Vol}(M)}+\alpha\frac{\int_M\langle\xi\cdot\Omega\cdot\psi,\psi\rangle v_g}{\mathrm{Vol}(M)}.
\ee
The Min-Max principle concludes the proof.

\hfill$\square$\\\\

\noindent{\bf Remark.} It would be simpler to compute $|D_M\psi|^2$ by (\ref{eqDMpsi}) to obtain the upper bound in Theorem \ref{lem:2}. In fact, we chose to compute $D_M^2\psi$ because its expression is required for the limiting case in the Min-Max principle.

\begin{eexemples}\label{exegalprodriem}
{\rm\noindent
\begin{enumerate}
\item Let $M^{n+1}$ be a compact spin manifold admitting a unit parallel vector field $\xi$. This is equivalent to the fact that $M$ is locally the product of two Riemannian manifolds, one of those being $1$-dimensional. The vector field $\xi$ defines a Riemannian flow on $M$ with totally geodesic fibres and vanishing O'Neill tensor. If moreover we assume that $M$ carries a $(0,\beta)$-transversal Killing spinor (for example $M:=\s^1\times N,$ where $N$ is a Riemannian spin manifold carrying a non-zero $\beta$-Killing spinor and $\s^1$ is endowed with the trivial spin structure), we deduce from \cite[Thm. 1.1]{AlexGrantchIvanov98} that $\lambda_1(D_M^2)= n^2\beta^2$, hence we have equality in \eqref{estimation}.

\item Consider an $n+1$-dimensional f\mbox{}lat torus $M=\mathbb{T}^{n+1}$ with its trivial spin structure. It is shown in \cite{GinHabib1} that a parallel spinor on $M$ induces an $(\alpha,0)$-transversal Killing spinor for some $\alpha\in\R^*$ (in fact $\alpha$ is related to the length of the $\s^1$-f\mbox{}ibers such that $\mathbb{T}^{n+1}=\s^1\times \mathbb{T}^n$). In this case, the upper bound given in Theorem \ref{lem:2} is equal to $\alpha^2$, which is positive. Therefore the equality is in that case not attained, since the first eigenvalue is equal to $0.$ This shows that the upper bound obtained in Theorem \ref{lem:2} depends sensitively on the transversal Killing spinor chosen on the Riemannian f\mbox{}low.
\end{enumerate}
}
\end{eexemples}

\section{Case of Sasakian manifolds}\label{sec:Sasaki}
\setcounter{equation}{0}
\noindent
In this section, we simplify the estimate in Theorem \ref{lem:2} for Sasakian manifolds and we study the limiting case. We illustrate our results on examples.

\subsection{The estimate}

\noindent In \cite{GinHabib1} it is shown that the existence on a Sasakian manifold $M^{2m+1}$ of a non-zero $(\alpha,\beta)$-transversal Killing spinor $\psi$ with real $\alpha$ and $\beta$ implies $\alpha\beta=0.$ If moreover $M$ is compact and $m>1,$ then it cannot carry any non-zero $(0,\beta)$-transversal Killing spinor for non-zero $\beta$. 

\begin{cor}\label{csasaki}
Under the hypotheses of {\rm Theorem \ref{lem:2}}, if furthermore $M$ is Sasakian and $\alpha,\beta\in\R$ then
\begin{equation}\label{eq:estimsasaki}
\lambda_1(D_M^2)\leq\alpha^2+4m^2\beta^2+\frac{m^2}{4}+\alpha\frac{\int_M\langle\xi\cdot\Omega\cdot\psi,\psi\rangle v_g}{\mathrm{Vol}(M)}.
\end{equation}
Moreover, the equality can only occur if $\beta=0$. In that case, $\psi=\psi_0+\psi_m$, the f\mbox{}low is $\eta$-Einstein and if moreover $\alpha\neq 0$ then either $\psi$ is eigen for the Clif\mbox{}ford action by $\Omega$ or $m$ is odd.
\end{cor}

\noindent{\bf Proof.}
From (\ref{eqdecSigmaSasaki}), the following estimate holds pointwise: $|\Omega\cdot\psi|^2\leq m^2|\psi|^2$. 
Together with (\ref{estimation}) this straightforward yields the inequality (\ref{eq:estimsasaki}).\\
If (\ref{eq:estimsasaki}) is an equality then on the one hand $|\Omega\cdot\psi|^2=m^2|\psi|^2$, from which $\psi=\psi_0+\psi_m$ follows, and on the other hand $D_M^2\psi=\lambda_1(D_M^2)\psi$, that is, using (\ref{eq:DM2psi}) with $\kappa=0$, $\nabla h=0$ and $\alpha\beta=0$,
\be
\lambda_1(D_M^2)\psi&=&(\alpha^2+4m^2\beta^2)\psi-\frac{1}{4}\Omega\cdot\Omega\cdot\psi+\alpha\xi\cdot\Omega\cdot\psi+2\beta\Omega\cdot\psi\\
&=&(\alpha^2+4m^2\beta^2+\frac{m^2}{4})\psi-m\alpha\psi_0+(-1)^m m\alpha\psi_m\\&&-2mi\beta(\psi_0-\psi_m),
\ee
therefore $\beta=0$. The rest of the proof directly follows from \cite[Prop. 5.1]{GinHabib1}.
\hfill$\square$\\

\noindent{\bf Remark.} Decomposing $\psi$ into eigenspinors for the Clif\mbox{}ford action by $\Omega$ it is actually easy to show that
\[\lambda_1(D_M^2)\leq 4m^2\beta^2+\max_{0\leq r\leq m}\Big(\alpha+(-1)^r\frac{(2r-m)}{2}\Big)^2,\]
which has the advantage to be independent of $\psi$. However, the equality case in that estimate is more complicated to handle as  (\ref{eq:estimsasaki}), that is why we shall not consider it further on in this paper.\\

\subsection{$\mathcal{D}$-homothetic changes of metric} 
We now show the existence of non-zero harmonic spinors under some sui\-ta\-ble conditions on a Sasakian manifold carrying an $(\alpha,0)$-transversal Killing spinor. Recall that a $\mathcal{D}$-homothetic deformation of $g$ on a Sasakian manifold $(M^{2m+1},g,\xi)$ is a metric of the form $g_t:=t^2g|_\xi+tg|_Q$ on $M$ for some real number $t>0$ \cite{Tan68,GinHabib1}. The manifold $(M,g_t,\xi_t:=\frac{1}{t}\xi)$ is again Sasakian and, if furthermore $(M,g)$ is spin and carries an $(\alpha,\beta)$-transversal Killing spinor, then so is $(M,g_t)$ and carries an $(\frac{\alpha}{t},\frac{\beta}{\sqrt{t}})$-transversal Killing spinor \cite[Lemma 3.3]{GinHabib1}. 

\begin{pro}\label{pexistharmspinors}
Let $(M^{2m+1},g,\xi)$ be a compact Sasakian spin manifold car\-ry\-ing a non-zero $(\alpha,0)$-transversal Killing spinor $\psi$ with $\alpha\neq 0.$ Assume that the equality case in {\rm Corollary \ref{csasaki}} is attained on $(M^{2m+1},g,\xi).$
\beit\item If $\psi_0\neq 0$ and $\alpha>0$ then the spinor f\mbox{}ield $\ovl{\psi_0}^t$ corresponding to $\psi_0$ on the Sasakian manifold $(M^{2m+1},\ovl{g}_t,\ovl{\xi}^t)$ is harmonic for $t=\frac{2\alpha}{m}$.
\item If $\psi_m\neq 0$ and $(-1)^{m+1}\alpha>0$ then the spinor f\mbox{}ield $\ovl{\psi_m}^t$ correspon\-ding to $\psi_m$ on the Sasakian manifold $(M^{2m+1},\ovl{g}_t,\ovl{\xi}^t)$ is harmonic for $t=(-1)^{m+1}\frac{2\alpha}{m}$.
\eeit
\end{pro}

\noindent{\bf Proof.} We already know from \cite[Lemma 3.3]{GinHabib1} that the spinor f\mbox{}ield $\ovl{\psi}^t$ is an $(\frac{\alpha}{t},0)$-transversal Killing spinor. 
It now follows from (\ref{eqDM2psinablah}) that, for the metric $\ovl{g}_t$,
\[\ovl{D}_t^2\ovl{\psi}^t=\frac{\alpha^2}{t^2}\ovl{\psi}^t-\frac{1}{4}\,\ovl{\Omega}^t\,\ovl{\cdot}\,\ovl{\Omega}^t\,\ovl{\cdot}\,\ovl{\psi}^t+\frac{\alpha}{t}\,\ovl{\xi}^t\ovl{\cdot}\,\ovl{\Omega}^t\,\ovl{\cdot}\,\ovl{\psi}^t,\]
where $\ovl{D}_t$ is the Dirac operator associated with $\ovl{g}_t.$ Using \cite[Lemma 2.3]{GinHabib1} we have $\ovl{\Omega}^t\,\ovl{\cdot}\,\ovl{\psi}^t=\ovl{(\Omega\cdot\psi)}^t$ and we obtain
\[ \ovl{D}_t^2\ovl{\psi}^t=\frac{\alpha^2}{t^2}\ovl{\psi}^t-\frac{1}{4}\ovl{(\Omega\cdot\Omega\cdot\psi)}^t+\frac{\alpha}{t}\,\ovl{(\xi\cdot\Omega\cdot\psi)}^t.\]
Assume now that $M$ is compact and (\ref{eq:estimsasaki}) is an equality on $(M^{2m+1},g,\xi)$. If $\alpha=0$ then $\psi=\psi_0+\psi_m$ and we get $\ovl{D}_t^2\ovl{\psi}^t=\frac{m^2}{4}\ovl{\psi}^t,$ which implies that $\ovl{\psi}^t$ is harmonic for no $t.$ If now $\alpha\neq 0$ we get for $\psi=\psi_0$ 
\be
\ovl{D}_t^2\ovl{\psi}^t&=&\frac{\alpha^2}{t^2}\ovl{\psi}^t+\frac{m^2}{4}\ovl{\psi}^t-m\frac{\alpha}{t}\ovl{\psi}^t\\
&=&(\frac{m}{2}-\frac{\alpha}{t})^2\ovl{\psi}^t,
\ee
which vanishes for $t=\frac{2\alpha}{m}$, which can of course only happen if $\alpha>0$. In the case $\psi=\psi_m$ the equality above becomes  
\be
\ovl{D}_t^2\ovl{\psi}^t&=&\frac{\alpha^2}{t^2}\ovl{\psi}^t+\frac{m^2}{4}\ovl{\psi}^t+ (-1)^m m\frac{\alpha}{t}\ovl{\psi}^t\\
&=&(\frac{m}{2}+(-1)^m\frac{\alpha}{t})^2\ovl{\psi}^t,
\ee
which vanishes for $t=(-1)^{m+1}\frac{2\alpha}{m}.$ This concludes the proof. 
\hfill$\square$\\

\begin{eexemples}\label{exemplesasaki}
{\rm\noindent\begin{enumerate}
 \item Let $(\s^{2m+1},g)$ be the round sphere equipped with the round metric $g$ with sectional curvature equal to $1.$ It is a Sasakian manifold w.r.t. the vector field $\xi_x=ix$ for every $x\in \s^{2m+1}\subset \C^{m+1}.$ Since the sphere admits non-zero $\pm\frac{1}{2}$-Killing spinors, the first eigenvalue of the square of the Dirac operator is equal to $\frac{(2m+1)^2}{4}$ \cite{Fried80}.
Moreover it is proved in \cite[Prop. 5.7]{GinHabib1} that, whatever the parity of $m$ is, there always exists a $1$-dimensional space of $(-\frac{m+1}{2},0)$-transversal Killing spinors lying in $\Gamma(\Sigma_0 M)$. In particular
\[(\frac{m}{2}-\alpha)^2=(\frac{m}{2}+\frac{m+1}{2})^2=\frac{(2m+1)^2}{4},\]
hence (\ref{eq:estimsasaki}) is an equality.   
\item For any $\mathcal{D}$-homothetic change $\ovl{g}_t$ of the standard metric $g$ on $\s^{2m+1}$ the number $(\frac{m}{2}+\frac{m+1}{2t})^2$ is an eigenvalue of the square of the Dirac operator on $(\s^{2m+1},\ovl{g}_t)$. However since there exists at least one $t\in\R_+^*$ such that $(\s^{2m+1},\ovl{g}_t)$ admits non-zero harmonic spinors \cite{Baer96} the equality in (\ref{eq:estimsasaki}) doesn't survive under arbitrary variations of $t$. The eigenvalue we obtain from the existence of transversal Killing spinors,  $(\frac{m}{2}+\frac{m+1}{2t})^2$, is not the one that crosses the zero line. In an equivalent way, Proposition \ref{pexistharmspinors} cannot be applied to the existence of harmonic spinors on $(\s^{2m+1},\ovl{g}_t)$ since the assumptions on the sign of $\alpha$ are not fulf\mbox{}illed.   
\item 
Consider now the quotient $M:=\Gamma\setminus \s^3$ where $\Gamma$ is a non-trivial f\mbox{}inite subgroup of $\SU_2$. 
Every such quotient endowed with the metric $g$ induced by the standard metric on $\s^3$ with constant sectional curvature $1$ is again a Sasakian manifold.
Moreover (see \cite[Notes 5.9.2]{GinHabib1}) there exists a spin structure on $M$ for which $M$ carries a one-dimensional space of $(-1,0)$-transversal Killing spinors and analogously the upper bound \eqref{eq:estimsasaki} is attained. 
\item Let 
$M_r:=\vspace{-1mm}{\Gamma_r}\hspace{-1mm}\setminus\hspace{-0.5mm}\vspace{1mm} G$ be a Heisenberg manifold as in \cite[Ex. 3.9.5]{GinHabib1}. 
It is a Sasakian ma\-ni\-fold carrying a $2$-dimensional space of transversal parallel spinors w.r.t. its canonical spin structure \cite[Sec. 3]{GinHabib1}. 
On the other hand, B. Ammann and C. B\"ar proved in \cite{AmmBaer98} that the smallest (in absolute value) eigenvalue of the Dirac o\-pe\-ra\-tor for any spin structure is equal to $\pm\frac{1}{2}$. Hence the limiting case in Corollary \ref{csasaki} is attained.
\end{enumerate}
}
\end{eexemples}

\section{The 3-dimensional case} \label{pro:dimentrois}
\setcounter{equation}{0}
\noindent
In this section, we consider $3$-dimensional Riemannian flows. We deduce from Theorem \ref{lem:2} an estimate in terms of the O'Neill tensor and we classify all the manifolds satisfying the equality case.

\subsection{The estimate}

\noindent Let $(M^3,g,\mathcal{F})$ be a $3$-dimensional Riemannian f\mbox{}low. Since the map $h$ is a skew-symmetric endomorphism-field of $Q$, it can be represented locally by the matrix
$$\begin{pmatrix}
0&-b\\b&0
\end{pmatrix},$$
where $b: M \longrightarrow \R$ is a smooth real-valued function defined on $M$. We recall that the complex volume form
$\omega_3=-\xi\cdot e_1\cdot e_2$ (where $\{\xi, e_1, e_2\}$ is a local oriented frame of $TM$) acts as the identity on the spinor bundle. Note that the existence of a non-zero $(\alpha,\beta)$-transversal Killing spinor for real $\alpha,\beta$ on a minimal compact flow $M$ implies that $\alpha\beta=0$ and $\xi(b)=0$ \cite[Prop. 6.2]{GinHabib1}.
The scalar curvature of $M$ is equal to
\begin{equation}
{\rm Scal}_M=2(4\beta^2-b^2-4\alpha b).
\label{eq:scal}
\end{equation} 
  
\begin{cor}\label{pestimdim3}
Under the hypotheses of {\rm Theorem \ref{lem:2}}, assume furthermore that $\alpha,\beta\in\R$ and $\kappa=0$. Then 
\begin{equation} \label{ineqlambda1dim3}
\lambda_1(D_M^2)\leq 4\beta^2+\frac{1}{\mathrm{Vol}(M)}\int_M(\frac{b}{2}-\alpha)^2 v_g.
\end{equation}
If moreover \eqref{ineqlambda1dim3} is an equality, then $b$ is constant and either $\beta=0$ or $b=0.$ 
\end{cor}

\noindent{\bf Proof.}
The estimate is a direct consequence of Theorem \ref{lem:2} and from the fact that for any spinor field $\psi\in \Gamma(\Sigma M)$ we have $\Omega\cdot\psi=b\xi\cdot\psi.$
Now if (\ref{ineqlambda1dim3}) is an equality, then we have $D_M^2\psi=\lambda_1(D_M^2)\psi$. On the other hand \eqref{eqDM2psinablah} gives
\be
D_M^2\psi&=&(\alpha^2+4\beta^2)\psi+2b\beta\xi\cdot\psi+\frac{1}{2}\nabla_\xi^M\Omega\cdot\psi\\
& &+\frac{1}{4}\sum_{j,k=1}^2 \xi\cdot e_j\cdot e_k\cdot\nabla_{e_j}h(e_k)\cdot\psi+\frac{b^2}{4}\psi-b\alpha\psi.
\ee
Using the fact that $b$ is constant along the flow, we then write
$$\nabla^M_\xi\Omega\cdot\psi=\nabla^M_\xi(b\xi)\cdot\psi=\xi(b)\xi\cdot\psi=0.$$
We also compute
\be
\sum_{j,k=1}^2\xi\cdot e_j\cdot e_k\cdot\nabla_{e_j}h(e_k)\cdot\psi&
=&-e_1(b)\xi\cdot e_2\cdot\psi-e_1(b)\xi\cdot e_2\cdot\psi\\
&&+e_2(b)\xi\cdot e_1\cdot\psi+e_2(b)\xi\cdot e_1\cdot\psi\\
&=&2db\cdot\psi.
\ee
Hence, we find
$$\lambda_1(D_M^2)\psi=(\alpha^2+4\beta^2+\frac{b^2}{4}-b\alpha)\psi+2b\beta\xi\cdot\psi+\frac{1}{2}db\cdot\psi.$$
The Hermitian product with $\psi$ gives after identifying the real parts that
\[\lambda_1(D_M^2)=4\beta^2+(\frac{b}{2}-\alpha)^2,\] 
which implies that the function $b$ must be constant on $M.$ Hence $b\beta=0$ and the proof is achieved.
\hfill$\square$

\begin{eexemples}
{\rm\noindent\begin{enumerate}
\item Let $M$ be the Riemannian product $\s^1\times \s^2$, where $\s^2$ carries its canonical spin structure and its canonical metric with sectional curvature $1$ and $\s^1$ its trivial spin structure. The manifold $M$ is then a trivial fibration over $\s^2$ with vanishing O'Neill tensor. Moreover, we have already shown in Example \ref{exegalprodriem} that the first eigenvalue of the Dirac operator is equal to $1$ and the manifold $M$ carries a $(0,\pm\frac{1}{2})$-transversal Killing spinor. Hence we have $$\lambda_1(D_M^2)=1=4\beta^2+\underbrace{(\frac{b}{2}-\alpha}_{0})^2,$$
and the upper bound in \eqref{ineqlambda1dim3} is attained.

\item Let $M=\Gamma\setminus \s^3$ as in Examples \ref{exemplesasaki}. Since $M$ is a Sasakian manifold one may choose $b=1.$ Moreover, the first eigenvalue of the square of the Dirac operator of $M$ is equal to $\frac{9}{4}$ and it admits a $(-1,0)$-transversal Killing spinor. Then we compute $$\lambda_1(D_M^2)=\frac{9}{4}=\underbrace{4\beta^2}_{0}+\underbrace{(\frac{b}{2}-\alpha}_{\frac{3}{2}})^2,$$
and we have equality in \eqref{ineqlambda1dim3}.

\item Let $M_r$ be a Heisenberg manifold (see Examples \ref{exemplesasaki}). The manifold $M_r$ is a Sasakian manifold carrying a $(0,0)$-transversal Killing spinor and the first eigenvalue is equal to $\frac{1}{4}.$ Therefore, we find 
$$\lambda_1(D_M^2)=\frac{1}{4}=\underbrace{4\beta^2}_{0}+\underbrace{(\frac{b}{2}-\alpha}_{\frac{1}{2}})^2,$$
which gives the equality in \eqref{ineqlambda1dim3}.


\end{enumerate}
}
\end{eexemples}

\subsection{Classification of limiting manifolds}

\noindent In this section, we classify all $3$-dimensional manifolds sa\-tis\-fying the limiting case in Inequality \eqref{ineqlambda1dim3}. Recall that, if we have equality, then either $b=0$ (which means that $M$ is a local a Riemannian product) or $\beta=0$ with non-zero constant $b$ (which implies that $(M,b^2g)$ is Sasakian). Using the classification of compact $3$-dimensional Sasakian manifolds \cite{Geiges97,Belgun00} we classify in \cite{GinHabib1}, up to $\mathcal{D}$-homothetic deformations of the metric, all compact manifolds admitting $(\alpha,\beta)$-transversal Killing spinors with constant $b.$ Therefore, we discuss the following cases \cite[Prop. 6.4]{GinHabib1}:
\begin{itemize}
\item {\bf Case where $\alpha=\beta=0$}: The manifold $M$ is either isometric to the flat torus $\T^3$ or to the Heisenberg manifold, both endowed with their trivial spin structure.

\item {\bf Case where $\alpha=0, \beta\neq 0$}: In this case, the O'Neill tensor vanishes and the manifold is isometric to the Riemannian product $\s^1\times \s^2$, where $\s^1$ carries its canonical spin structure and $\s^2$ its canonical one.

\item {\bf Case where $\alpha\neq 0, \beta=0$}: Since $b$ is constant, we distinguish two cases:
\begin{itemize}
\item First, if $b\neq 0$ then the scalar curvature on $(M,b^2g)$ is equal by \eqref{eq:scal} to $-2-8\alpha$. The manifold $M$ is then either isometric to the quotient $\Gamma\setminus \s^3$ ($\alpha<0$) or the manifold $\Gamma\setminus \widetilde{\rm{PSL}_2}(\R)$ ($\alpha>0$). 
\item If $b$ is equal to zero, then the manifold $M$ is Ricci flat and in dimension $3$ this implies flatness. Hence $M$ is a Bieberbach manifold with a suitable spin structure, i.e., $M$ is isometric to the quotient $G_i\setminus \R^3, i=1,\ldots,5$, where $G_i$ are subgroups of the group of Euclidean motions ${\rm SO}_3\ltimes \R^3$.
\end{itemize}  
Note that in \cite{GinHabib1}, we show that an $(\alpha,0)$-transversal Killing spinor on a Bieberbach manifold $M$ could be lifted to an $(\alpha,0)$-transversal Killing spinor on the Euclidian space $\R^3,$ the universal cover on $M,$ invariant under the action of $G_i.$ Therefore the existence of such spinors induces equivariance conditions on the generators of $G_i,$ for $i=1,\ldots,5$ (they are expressed below). Now we recall these conditions and we will show that the first eigenvalue of the square of the Dirac operator on these manifolds is equal to $\alpha^2.$\\

\noindent The spin structures on $G_i\setminus \R^3, i=1,\ldots,5,$ are in bijective relation to the homomorphisms $\varepsilon:G_i\rightarrow \spin_3$ such that $r=\ad\circ\varepsilon,$ i.e., the following diagram commutes \cite{Pfaeffle00}

$$
\begin{diagram}
\node{}\node{\spin_3}\arrow{s,r}{\ad}\\
\node{G_i}\arrow{ne,t}{\varepsilon}\arrow{e,t}{r}\node{{\rm SO}_3}
\end{diagram}
$$

where $r$ is defined for all $g=(A,a)\in G_i$ by $r(A,a)=A.$ Next we denote by $(\delta_1,\delta_2,\delta_3)\in \{0,1\}^3$ and by $\{e_1,e_2,e_3\}$ the canonical basis in $\R^3$ and let $\tilde\xi$ be the lift of $\xi$ to $\R^3.$\\ 
 
{\bf $\bullet$ Case  $i=1$}: The group $G_1$ is generated by three translations associated with the basis $\{a_1,a_2,a_3\}$ of $\R^3.$ The values of $\varepsilon$ on the ge\-ne\-ra\-tors are given by $\varepsilon(a_j)=e^{i\pi\delta_j},$ for $j=1,2,3.$ The manifold $M$ is then a torus carrying $8$ spin structures. The existence of an $(\alpha,0)$-transversal Killing spinor on $M$ implies $\langle a_j,\tilde\xi\rangle\equiv \frac{\pi}{\alpha}\delta_j\; [\frac{2\pi}{\alpha}].$ For the spin structure $(\delta_1,\delta_2,\delta_3)=(1,0,0)$ and for $\tilde\xi=e_1$ consider the basis $a_1=(\frac{\pi}{\alpha},0,0),\,a_2=(0,1,0),\,a_3=(0,0,1).$
Therefore, associated with this basis, the manifold $M$ carries an $(\alpha,0)$-transversal Killing spinor and moreover the spectrum of the square of the Dirac operator on $M$ is given by \cite{Fried84}
$$\{4\alpha^2(m+\frac{1}{2})^2+4\pi^2n^2+4\pi^2p^2|\,\, (m,n,p)\in \Z^3\}.$$
Hence the lowest eigenvalue for the square of the Dirac operator is associated with $(m,n,p)=(0,0,0)$ and it is equal to $\alpha^2$ which is the upper bound in \eqref{ineqlambda1dim3}. Finally, we point out that for the trivial spin structure $(\delta_1,\delta_2,\delta_3)=(0,0,0)$, the Dirac operator admits $0$ as an eigenvalue and the limiting case in \eqref{ineqlambda1dim3} could not be achieved.\\

{\bf $\bullet$ Case $i=2$}: The group $G_2$ is generated by three translations associated with the basis  
$a_1=(0,0,H),\, a_2=(L,0,0),\, a_3=(T,S,0)$ with $H,L,S>0,\,T\in \R,$
and by $g=(A,\frac{1}{2}a_1)$ where $A$ is a $\pi$-rotation around $z$-axis. The values of $\varepsilon$ on the generators are given by $\varepsilon(a_1)=-1,\,\,\varepsilon(a_j)=e^{i\pi\delta_j}$ for $j=2,3,$ and $\varepsilon(g)=e^{i\pi\delta_1} e_1\cdot e_2.$ The existence of an $(\alpha,0)$-transversal Killing spinor on $M$ implies that, $\delta_j=0$ for $j=2,3,$ and $H\equiv\frac{\pi}{\alpha}(1+2\delta_1)\; [\frac{4\pi}{\alpha}].$ The spectrum of the square of the Dirac operator is given for the spin structure $(\delta_1,\delta_2,\delta_3)=(0,0,0)$ by 
\begin{equation*}
\{\frac{4\pi^2}{H^2}(k+\frac{1}{2})^2+\frac{4\pi^2l^2}{L^2}+\frac{4\pi^2}{S^2}(m-\frac{T}{L}l)^2,\,\,\frac{4\pi^2}{H^2}(2\mu+\frac{1}{2})^2\},
\end{equation*}
and for the spin structure $(\delta_1,\delta_2,\delta_3)=(1,0,0)$ by
\begin{equation*}
\{\frac{4\pi^2}{H^2}(k+\frac{1}{2})^2+\frac{4\pi^2l^2}{L^2}+\frac{4\pi^2}{S^2}(m-\frac{T}{L}l)^2,\,\,\frac{4\pi^2}{H^2}(2\mu+\frac{3}{2})^2\},
\end{equation*}
where $\mu\in \Z$ and $(k,l,m)\in I$ with
$$I=\{(k,l,m)\in \Z^3|\,\, m\geq 1\}\cup \{(k,l,0)\in \Z^3|\,\, l\geq 1\}.$$ 
Hence for $\delta_1=0$ (resp. $\delta_1=1$), the lowest eigenvalue corresponds to $k=l=\mu=0$ (resp. $k=l=0$) and $m=1$ (resp. $m=1$ and $\mu=-1$) and it is equal to $\frac{\pi^2}{H^2}.$ We then deduce that for $H=\frac{\pi}{\alpha}$ (resp. $H=-\frac{\pi}{\alpha}$), the equality in \eqref{ineqlambda1dim3} is realized.\\

{\bf $\bullet$ Case $i=3$}: The group $G_3$ is generated by three translations associated with the basis
$a_1=(0,0,H),\, a_2=(L,0,0),\, a_3=(-\frac{1}{2}L,\frac{\sqrt{3}}{2}L,0)$ with $H,L>0,$
and by $g=(A,\frac{1}{3}a_1)$ where $A$ is a $\frac{2\pi}{3}$-rotation around $z$-axis. The values of $\varepsilon$ on the generators are given by $\varepsilon(a_1)=-e^{i\pi\delta_1}$, $\varepsilon(a_j)=1$ for $j=2,3,$ and $\varepsilon(g)=e^{i\pi\delta_1}(\frac{1}{2}+\frac{\sqrt{3}}{2} e_1\cdot e_2).$ 
The existence of an $(\alpha,0)$-transversal Killing spinor implies that $H\equiv\frac{\pi}{\alpha}(1+3\delta_1) \;[\frac{6\pi}{\alpha}]$. The spectrum of the square of the Dirac operator is given for $\delta_1=0,$ by 
$$\{\frac{4\pi^2}{H^2}(k+\frac{1}{2})^2+\frac{4\pi^2l^2}{L^2}+\frac{4\pi^2}{3L^2}(l-2m)^2,\,\,\frac{4\pi^2}{H^2}(4\mu+\frac{1}{2})^2\},$$
and for $\delta_1=1,$ by 
$$\{\frac{4\pi^2k^2}{H^2}+\frac{4\pi^2l^2}{L^2}+\frac{4\pi^2}{3L^2}(l-2m)^2,\,\,\frac{4\pi^2}{H^2}(3\mu+2)^2\},$$
where $\mu\in \Z$ and $(k,l,m)\in I$ with
$$I=\{(k,l,m)\in \Z^3|\,\, l\geq 1, m=0,\ldots,l-1\}.$$ 
Hence for $\delta_1=0,$ the lowest eigenvalue is equal to $\frac{\pi^2}{H^2}$ which implies the equality in \eqref{ineqlambda1dim3} for $H=\frac{\pi}{\alpha}.$ Now for $\delta_1=1,$ the lowest eigenvalue is equal to ${\rm min}(\frac{16\pi^2}{3L^2},\frac{4\pi^2}{H^2}).$ Thus if $H^2>\frac{3}{4}L^2,$ the minimum is then equal to $\frac{4\pi^2}{H^2}$ and the equality is attained for $H=-\frac{2\pi}{\alpha}.$\\

{\bf $\bullet$ Case $i=4$}: The group $G_4$ is generated by three translations associated with the basis  
$a_1=(0,0,H),\, a_2=(L,0,0),\, a_3=(0,L,0)$ with $H,L>0,$
and by $g=(A,\frac{1}{4}a_1),$ where $A$ is a $\frac{\pi}{2}$-rotation around $z$-axis. The values of $\varepsilon$ on the generators are given by $\varepsilon(a_1)=-1,\,\,\varepsilon(a_j)=e^{i\pi\delta_2}$ for $j=2,3,$ and $\varepsilon(g)=e^{i\pi\delta_1}(\frac{\sqrt{2}}{2}+\frac{\sqrt{2}}{2} e_1\cdot e_2).$ 
The existence of an $(\alpha,0)$-transversal Killing spinor implies that $\delta_2=0$ and $H\equiv \frac{\pi}{\alpha}(1+4\delta_1)\; [\frac{8\pi}{\alpha}]$. The spectrum of the square of the Dirac operator is given for $(\delta_1,\delta_2)=(0,0)$ by 
$$\{\frac{4\pi^2}{H^2}(k+\frac{1}{2})^2+\frac{4\pi^2}{L^2}(l^2+(m-l)^2),\,\,\frac{4\pi^2}{H^2}(4\mu+\frac{1}{2})^2\},$$
and for the spin structure $(\delta_1,\delta_2)=(1,0)$ by
$$\{\frac{4\pi^2}{H^2}(k+\frac{1}{2})^2+\frac{4\pi^2}{L^2}(l^2+(m-l)^2),\,\,\frac{4\pi^2}{H^2}(4\mu+\frac{5}{2})^2\},$$
where $\mu\in \Z$ and $(k,l,m)\in I$ with
$$I=\{(k,l,m)\in \Z^3|\,\, l\geq 1, m=0,\ldots,2l-1\}.$$ 
The lowest eigenvalue is given by $\frac{\pi^2}{H^2}$ for $\delta_1=0$ and we have equality in \eqref{ineqlambda1dim3} for $H=\frac{\pi}{\alpha}.$ Also for $\delta_1=1$ the lowest eigenvalue is equal to ${\rm min}(\frac{\pi^2}{H^2}+\frac{4\pi^2}{L^2},\frac{9\pi^2}{H^2}).$ Hence if $H^2>2L^2,$ the equality is achieved for $H=-\frac{3\pi}{\alpha}.$\\

{\bf $\bullet$ Case $i=5$}: The group $G_5$ is generated by three translations associated with the basis 
$a_1=(0,0,H),\, a_2=(L,0,0),\, a_3=(\frac{1}{2}L,\frac{\sqrt{3}}{2}L,0)$ with $H,L>0,$
and by $g=(A,\frac{1}{6}a_1)$ where $A$ is a $\frac{\pi}{3}$-rotation around the $z$-axis. The values of $\varepsilon$ on the generators are given by $\varepsilon(a_1)=-1,\,\,\varepsilon(a_j)=1$ for $j=2,3,$ and $\varepsilon(g)=e^{i\pi\delta_1}(\frac{\sqrt{3}}{2}+\frac{1}{2} e_1\cdot e_2).$ 
The existence of an $(\alpha,0)$-transversal Killing spinor implies that $H\equiv \frac{\pi}{\alpha}(1+6\delta_1) \;[\frac{12\pi}{\alpha}].$ The spectrum of the square of the Dirac operator is given for $\delta_1=0$ by 
$$\{\frac{4\pi^2}{H^2}(k+\frac{1}{2})^2+\frac{4\pi^2l^2}{L^2}+\frac{4\pi^2}{3L^2}(2l-m)^2,\,\,\frac{4\pi^2}{H^2}(6\mu+\frac{1}{2})^2\},$$
and for $\delta_1=1$ by
$$\{\frac{4\pi^2}{H^2}(k+\frac{1}{2})^2+\frac{4\pi^2 l^2}{L^2}+\frac{4\pi^2}{3L^2}(2l-m)^2,\,\,\frac{4\pi^2}{H^2}(6\mu+\frac{7}{2})^2\},$$
where $\mu\in \Z$ and $(k,l,m)\in I$ with
$$I=\{(k,l,m)\in \Z^3|\,\, l\geq 1, m=0,\ldots,l-1\}.$$ 
The lowest eigenvalue is given by $\frac{\pi^2}{H^2}$ for $\delta_1=0$ and the equality in \eqref{ineqlambda1dim3} is attained for $H=\frac{\pi}{\alpha}.$ Also for $\delta_1=1$ the lowest eigenvalue is equal to ${\rm min}(\frac{25\pi^2}{H^2},\frac{\pi^2}{H^2}+\frac{28\pi^2}{3L^2}).$ Therefore, if $H^2>\frac{18}{7}L^2,$ the equality is realized for $H=-\frac{5\pi}{\alpha}.$
\hfill$\square$
\end{itemize}
To sum up, we prove that on each of the Bieberbach manifold $G_i\setminus\R^3$ for $i=1,\ldots,5,$ the upper bound in the estimate \eqref{ineqlambda1dim3} is attained.

\section{Comparison with a lower bound of the spectrum} \label{sec:comp}
\setcounter{equation}{0}

\noindent In this section we compare the upper bound obtained in Corollaries \ref{csasaki} and \ref{pestimdim3} with a lower bound of the spectrum of the Dirac operator analogous to that proved by O. Hijazi in \cite{Hij95} and generalized by the second author in \cite{Habibthese}. We recall first that if the scalar curvature of $(M^{n+1},g)$ is {\it negative} at a point (which is the case when e.g. $M$ is Sasakian and carries a transversally parallel spinor \cite{GinHabib1}), T. Friedrich's estimate \cite{Fried80} 
\[\lambda_1(D_M^2)\geq\frac{n+1}{4n}\inf_M{\mathrm{Scal}_M},\]
doesn't give any information on the spectrum of the Dirac operator. In \cite{Hij95} (the result was generalized by the second author in \cite{Habibthese}), O. Hijazi improves Friedrich's inequality by giving a new estimate depending on the eigenspinor but which may be non-trivial even if the scalar curvature of $M$ is negative at a point:

\begin{thm}[\cite{Hij95},\cite{Habibthese}]
Let $(M,g)$ be a compact spin manifold and $\phi$ be a non-zero eigenspinor of $D_M^2$ associated with the eigenvalue $\lambda_1(D_M^2).$ Then we have  
\begin{equation}
\lambda_1(D_M^2)\geq \inf_{M\setminus Z_\phi}(\frac{{\rm Scal}_M}{4}+|E^\phi|^2),
\label{eq:0034}
\end{equation}
where $Z_\phi:=\{x\!\in\! M\,|\,\phi_x\!=\!0\}$ and $E^\phi$ is the $2$-tensor f\mbox{}ield def\mbox{}ined on $M\setminus Z_\phi$ by
\[E^\phi(X,Y):=\Re(\langle Y\!\cdot\! \nabla^M_X\phi,\frac{\phi}{|\phi|^2}\rangle)\]
for all $X,Y\in \Gamma(TM)$.
\end{thm}

\noindent The proof of (\ref{eq:0034}) is based on the modif\mbox{}ication of the Levi-Civita connection in the direction of $E^\phi$ and the use of the Schr\"odinger-Lichnerowicz formula. The limiting case of (\ref{eq:0034}) is characterized by the spinor f\mbox{}ield $\phi$ satisfying for all $X\in \Gamma(TM)$ the equation
\begin{equation}\label{eqEMpsi}
\nabla^M_X\phi=-E^\phi(X)\cdot\phi.
\end{equation}
In this case, the manifold doesn't carry any particular geometry, since the lower bound depends on the eigenspinor.\\

\noindent
Now, we will prove that if the limiting case in Inequalities \eqref{eq:estimsasaki} and \eqref{ineqlambda1dim3} is achieved, then so is in \eqref{eq:0034}. Let us first consider the $3$-dimensional case.
 
\begin{pro} \label{estimdime3} Let $(M,g,\mathcal{F})$ be a $3$-dimensional compact Riemannian spin flow carrying a non-zero $(\alpha,\beta)$-transversal Killing spinor $\psi$ for real $\alpha$ and $\beta$. Assume the first eigenvalue of the square of the Dirac operator to be equal to $4\beta^2+(\frac{b}{2}-\alpha)^2.$ Then the equality in \eqref{eq:0034} is always attained.
\end{pro}
{\bf Proof.} Using the fact that the complex volume form acts as the identity,  we have from Equations \eqref{eq:25} on $M$
$$
\left\{\begin{array}{ll}
\nabla^M_\xi\psi=(\alpha+\frac{b}{2})\xi\cdot\psi,\\\\ 
\nabla^M_{e_1}\psi=\beta e_2\cdot\psi-\frac{b}{2}e_1\cdot\psi,\\\\ 
\nabla^M_{e_2}\psi=-\beta e_1\cdot\psi-\frac{b}{2}e_2\cdot\psi,
\end{array}\right.
$$
where $\{\xi,e_1,e_2\}$ is a local basis of $TM.$ We easily compute the tensor $E^\psi$ and we find that $E^\psi=-(\frac{b}{2}+\alpha)\xi^\flat\otimes\xi+\frac{b}{2}{\rm Id}_Q-\beta J$, where $J$ is the canonical K\"ahler structure associated to the flow, see e.g. \cite[Sec. 6.1]{GinHabib1}. Hence the limiting case in \eqref{eq:0034} is achieved.
\hfill$\square$\\

\noindent{\bf Remarks.}
\noindent\ben\item Using the classification of manifolds in dimension $3$ established in Section \ref{pro:dimentrois}, we deduce that the following manifolds $G_i\setminus\R^3\, (i=1,\ldots,5),\break \Gamma\setminus \widetilde{\rm{PSL}_2}(\R)$ and the Heisenberg group are examples of $3$-dimensional manifolds satisfying the limiting case in \eqref{eq:0034} and in which the Friedrich estimate does not give any information on their Dirac spectrum, since they have either vanishing or negative scalar curvature.
\item Under the asumptions of Proposition \ref{estimdime3}, if $\alpha=0$, i.e., $M=\s^2\times\s^1,$ then the tensor $E^\psi$ is a skew-symmetric tensor since $E^\psi=-\beta J.$ 
\een

\noindent
Now, we will compare the estimate \eqref{eq:0034} to the one in Corollary \ref{csasaki}. We have the following theorem: 
\begin{pro}\label{tcomparHijazi}
Let $(M^{2m+1},g)$ be a compact Sasakian spin manifold carrying an $(\alpha,0)$-transversal Killing spinor $\psi.$ We assume that the limiting case of Inequality \eqref{eq:estimsasaki} is attained, i.e., that
$$\lambda_1(D_M^2)=\alpha^2+\frac{m^2}{4}+\alpha\frac{\int_M\langle\xi\cdot\Omega\cdot\psi,\psi\rangle v_g}{\mathrm{Vol}(M)}.$$
Then we have the following results:
\beit\item[i)] If $\alpha\neq 0$, the equality in \eqref{eq:0034} is always attained.
\item[ii)] If $\alpha=0,$ the equality in \eqref{eq:0034} is attained if and only if either $m$ is odd or $\psi\in\Gamma(\Sigma_0 M)\cup\Gamma(\Sigma_m M).$
\eeit
\end{pro}

\noindent{\bf Proof.} 
We already know that $D_M^2\psi=\lambda_1(D_M^2)\psi$. Moreover, the spinor $\psi$ is an $(\alpha,0)$-transversal Killing spinor lying pointwise in $\Sigma_0 M\oplus\Sigma_m M$, hence we have from (\ref{eqkilling}) and (\ref{eq:25})  
\begin{eqnarray}\label{eqnablaxiMpsi}
\nonumber \nabla_\xi^M\psi&=&\alpha\xi\cdot\psi+\frac{1}{2}\Omega\cdot\psi\\
\nonumber &=&\alpha\xi\cdot\psi+\frac{im}{2}(-\psi_0+\psi_m)\\
&=&\alpha\xi\cdot\psi+\frac{m}{2}\xi\cdot(\psi_0+(-1)^{m+1}\psi_m)
\end{eqnarray}
and, for every $Z\in\Gamma(Q)$,
\begin{eqnarray}\label{eqnablaZMpsi}
\nonumber \nabla_Z^M\psi&=&\frac{1}{2}\xi\cdot h(Z)\cdot\psi\\
\nonumber &=&\frac{i}{2}\xi\cdot Z\cdot(\psi_0-\psi_m)\\
&=&-\frac{Z}{2}\cdot(\psi_0+(-1)^{m+1}\psi_m).
\end{eqnarray}
We handle the proof of $i)$ and $ii)$ in three cases which contain all possibilities, see \cite[Prop. 5.1]{GinHabib1}.\\
{\bf $\bullet$ Case $m$ odd}: we deduce from (\ref{eqnablaxiMpsi}) and (\ref{eqnablaZMpsi}) that, for every $Z\in\Gamma(Q)$,
\[\left|\begin{array}{ll}
\nabla_\xi^M\psi&=(\frac{m}{2}+\alpha)\xi\cdot\psi\\
&\\
\nabla_Z^M\psi&=-\frac{Z}{2}\cdot\psi,
\end{array}\right.\]
so that $\psi$ satisf\mbox{}ies (\ref{eqEMpsi}) with $E^\psi=-(\alpha+\frac{m}{2})\xi^\flat\otimes\xi\oplus \frac{1}{2}\mathrm{Id}_Q$. Hence (\ref{eq:0034}) is an equality in that case.\\
{\bf $\bullet$ Case $\psi_0=0$ or $\psi_m=0$}: if e.g. $\psi_m=0$ (the subcase $\psi_0=0$ is analogous), then we deduce from (\ref{eqnablaxiMpsi}) and (\ref{eqnablaZMpsi}) exactly the same equations as in the f\mbox{}irst case and therefore the same result.\\
{\bf $\bullet$ Case $m$ even and $\alpha=0$}: using the local expression of $\Omega$ it can be proved that 
\[
\nabla_\xi^M\Omega\cdot\psi=h(\kappa)\cdot\xi\cdot\psi+\frac{1}{2}\sum_{j=1}^n e_j\cdot (\nabla_\xi h)(e_j)\cdot\psi,\]
in particular $\nabla_\xi^M\Omega=0$ on Sasakian manifolds, therefore (\ref{eqnablaxiMpsi}) can be rewritten under the form
\[\left|\begin{array}{ll}\nabla_\xi^M\psi_0&=\frac{m}{2}\xi\cdot\psi_0\\
&\\
\nabla_\xi^M\psi_m&=-\frac{m}{2}\xi\cdot\psi_m,
\end{array}\right.\]
so that $\psi$ is a solution of (\ref{eqEMpsi}) only if $\psi_0=0$ or $\psi_m=0$. If this is satisf\mbox{}ied, then the second case shows that (\ref{eq:0034}) is an equality. This concludes the proof.
\hfill$\square$

$ $\\

\noindent {\bf Remarks.}
{\rm\ben\item Under the assumptions of Proposition \ref{tcomparHijazi}, the estimate (\ref{eq:0034}) is in fact an equality for $\phi=\psi$ only if $\psi$ is already eigen for $D_M$. Furthermore, the tensor $E^\psi$ is under the same assumptions always symmetric.
\item Equation (\ref{eqEMpsi}) has been studied by C. B\"ar, P. Gauduchon and A. Moroianu who show that, if the tensor $E^\psi$ is symmetric and satisfies the Codazzi identity, then any solution of (\ref{eqEMpsi}) comes from a parallel spinor on the so-called generalized cylinder of the underlying manifold, see \cite[Thm. 8.1]{BaerGaudMor}. Their result may unfortunately not be applied to our transversal Killing spinors, since under the assumptions of Propositions \ref{estimdime3} and \ref{tcomparHijazi} the tensor $E^\psi$ may be neither symmetric nor of Codazzi-type. 
\een
}

\begin{eexemples}\label{excaslim}
{\rm\noindent\ben\item Let $M:=\s^3/\Gamma$ as in \cite[Notes 5.9]{GinHabib1} (we also allow $\Gamma=\{\mathrm{Id}\}$). We have already shown in Examples \ref{exemplesasaki}.1 and \ref{exemplesasaki}.3 that, for the standard metric with constant sectional curvature $1$ and some spin structure on $M$, the inequality (\ref{eq:estimsasaki}) is an equality. Since in that case $\alpha=-1$ we deduce from Proposition \ref{tcomparHijazi} $i)$ that (\ref{eq:0034}) is also an equality. This is not a surprise since in that case T. Friedrich's lower bound is already attained.
\item Let $(M_r,g,\xi)$ be a Heisenberg manifold as in \cite[Ex. 3.9.5]{GinHabib1}. We have already shown in Examples \ref{exemplesasaki}.4 that, for any of the standard metrics and any spin structure on $M$, the inequality (\ref{eq:estimsasaki}) is an equality. Since $m=1$ we deduce from Proposition \ref{tcomparHijazi} $ii)$ that (\ref{eq:0034}) is also an equality. Moreover equality in T. Friedrich's inequality cannot occur since $\mathrm{Scal}_{M_r}=-2$. This is to the knowledge of the authors the f\mbox{}irst example of compact spin manifold where (\ref{eq:0034}) is sharp but not T. Friedrich's inequality.
\een
}
\end{eexemples}

\bibliographystyle{amsplain}


\providecommand{\bysame}{\leavevmode\hbox to3em{\hrulefill}\thinspace}

\end{document}